\documentclass[a4paper]{amsart}

\usepackage[utf8]{inputenc}
\usepackage{amsrefs}
\usepackage{enumerate}
\usepackage{booktabs}
\usepackage{multirow}
\usepackage{color}
\usepackage{calc}

\newcommand{\set}[1]{\{\,#1\,\}}
\newcommand{\stirling}{\genfrac\{\}{0pt}{}}

\newtheorem{lemma}{Lemma}
\newtheorem{proposition}{Proposition}
\theoremstyle{definition}
\newtheorem*{definition}{Definition}

\renewcommand{\MR}[1]{}
\author{David Bremner}
\address{Faculty of Computer Science, University of New Brunswick\\
  Box 4400, Fredericton NB, E3B 5A3}
\email{bremner@unb.ca}
\author{Lars Schewe}
\address{Fachbereich Mathematik, TU Darmstadt\\ Schlossgartenstrasse
  7\\ 64289 Darmstadt}
\email{schewe@mathematik.tu-darmstadt.de}
\thanks{This research was carried out while the authors were guests of the 
Hausdorff Research Institute for Mathematics, Bonn.\\
The research of the first author was partly supported by NSERC Canada\\
Computational resources were supplied by the ACEnet consortium.}

\title{Edge-Graph Diameter Bounds for Convex Polytopes with Few Facets}
\begin{document}
\begin{abstract}
  We show that the edge graph of a $6$-dimensional polytope with $12$
  facets has diameter at most $6$, thus verifying the $d$-step
  conjecture of Klee and Walkup in the case of $d=6$. This implies
  that for all pairs $(d,n)$ with $n-d\leq 6$ the diameter of the edge
  graph of a $d$-polytope with $n$ facets is bounded by $6$, which
  proves the Hirsch conjecture for all $n-d\leq 6$. We show this
  result by showing this bound for a more general structure -- so-called
  matroid polytopes -- by reduction to a small number of
  satisfiability problems. 
\end{abstract}
\maketitle

It is a long-standing open problem to determine the maximal diameter
$\Delta(d,n)$ of a $d$-dimensional polytope with $n$ facets. Not much
is known even in small dimensions. The \emph{Hirsch Conjecture} states
that $\Delta(d,n)\leq n-d$. The special case that $\Delta(d,2d)\leq d$
is known as the \emph{$d$-step conjecture}. Showing the $d$-step
conjecture for a fixed $d$ implies (by an argument of Klee and
Walkup~\cite{MR0206823}) that the Hirsch conjecture holds for all
pairs $(n',d')$ with $n'-d'=d$. So establishing the $d$-step
conjecture for all $d$ implies the general Hirsch conjecture. So far,
the $d$-step conjecture has been shown for all $d\leq 5$ by Klee and
Walkup.

We show that the $d$-step conjecture is true in dimension $6$ as
well. We derive this result by considering a more general class of
objects, namely \emph{matroid polytopes}, i.e. oriented matroids,
which, if realizable, correspond to convex polytopes. We show that no
$6$-dimensional matroid polytope with $12$ vertices and a facet path
of length~$7$ exists. Then $\Delta(6,12)=6$ follows by considering
polarity and the already known bounds. 

To show that $\Delta(6,12)\leq 6$ we first give combinatorial
conditions for matroid polytopes that violate this bound. This
achieved through the study of path complexes (see Bremner et
al.~\cite{bremner.ea:pathcomplex}, cf. Section~\ref{sec:pathcomp}). We
then show that these conditions cannot be satisfied by an oriented
matroid. To show this we use a satisfiability solver to produce the
desired contradiction (see Section~\ref{sec:omsat}). We will use the
same method to show that $\Delta(4,11)=6$, which settles another
special case of the Hirsch conjecture. The latter result allows us to
also improve the upper bound on $\Delta(5,12)$ from 9 to 8. 

For small parameters there are already known general bounds that allow
us to compute or at least bound the diameter of polytopes. We
summarize them in the following Lemma. 
For an overview about the known bounds we refer to the respective
chapters in the books by Grünbaum~\cite{MR1976856} and
Ziegler~\cite{MR1311028} and the survey of Klee and
Kleinschmidt~\cite{MR913867}.

\definecolor{follow}{gray}{0.3}
\newcommand{\foll}[1]{\textcolor{follow}{#1}}
\newlength{\tmpbox}\setlength{\tmpbox}{\widthof{$[7,10]$}}
\newlength{\tmpboxfour}\setlength{\tmpboxfour}{(\tmpbox-\widthof{4})/2}
\newlength{\tmpboxfive}\setlength{\tmpboxfive}{(\tmpbox-\widthof{5})/2}
\newlength{\tmpboxsix}\setlength{\tmpboxsix}{(\tmpbox-\widthof{\{6,7\}})/2}
\begin{table}[t]
  \centering
  \caption{Bounds on $\Delta(d,n)$ known before the work described in this article.}
  \label{tab:known}
  \begin{tabular}{ll@{\qquad}c@{\,}c@{\,}c@{\,}c}
     &&\multicolumn{4}{c}{$n-d$}\\
&    & 4 & 5 & 6 & 7\\
\addlinespace
\multirow{4}{*}{$d$}& 4 & \hspace{\tmpboxfour}4\hspace{\tmpboxfour} & \hspace{\tmpboxfive}5\hspace{\tmpboxfive}  & \hspace{\tmpboxfive}5\hspace{\tmpboxfive}  & \hspace{\tmpboxsix}\{6,7\}\hspace{\tmpboxsix} \\
&5 & \foll{4} & 5 & 6 & [7,9] \\
&6 & \foll{4} & \foll{5} & \{6,7\} & [7,9]\\
&7 & \foll{4} & \foll{5} & \foll{\{6,7\}} & $[7,10]$\\
&
  \end{tabular}
\end{table}

\begin{lemma}[Klee~\cite{MR0165514}, Klee and Walkup~\cite{MR0206823},
  Holt~\cite{MR2074846}]
The following relations hold for the maximal diameter $\Delta(d,n)$ of
a $d$-polytope with $n$ facets: 
  \label{thm:simple}
  \begin{enumerate}
  \item $\Delta(3,n)=\left\lfloor \frac{2}{3}n\right\rfloor - 1$
  \item $\Delta(d,2d+k)\leq \Delta(d-1,2d+k-1)+\lfloor \frac{k}{2}\rfloor +
    1$ for all $d$ and $k=0,1,2,3$
  \item $\Delta(d,n)\leq \Delta(n-d,2(n-d))$ for all $(d,n)$
  \item $\Delta(d,n)=\Delta(n-d,2(n-d))$ for all $(d,n)$ with $n\leq
    2d$ \label{item:dstep}
  \item $\Delta(d,n)\geq n-d$ for all $n>d\geq 7$
  \end{enumerate}
\end{lemma}

Apart from these general results, some special cases have been solved
by Goodey~\cite{MR0310768}:
\begin{lemma}[Goodey~\cite{MR0310768}]
\label{thm:goodey}
The following relations hold for the maximal diameter $\Delta(d,n)$ of
a $d$-polytope with $n$ facets: 
  \begin{enumerate}
  \item $\Delta(4,10)=5$ and $\Delta(5,11)=6$
  \item $\Delta(6,13)\leq 9$ and $\Delta(7,14)\leq 10$
  \end{enumerate}
\end{lemma}

Table~\ref{tab:known} summarizes the bounds on $\Delta(d,n)$ that
follow from Lemmas~\ref{thm:simple}~and~\ref{thm:goodey}.  The values
printed in grey follow from the observation (\ref{item:dstep}) of
Lemma~\ref{thm:simple}.

It is not difficult to see (e.g. by a perturbation argument) that
$\Delta(d,n)$ is always attained by a simple polytope. Thus, it is
sufficient for our purposes to restrict our attention to this class of
polytopes.  It will also be useful to investigate the problem in a
polar setting. Thus, we will be looking at $d$-dimensional simplicial
polytopes with $n$ vertices. In this setting $\Delta(d,n)$ is just the
maximal length of a shortest facet path in the polytope.

The rest of this article is organized as follows: we will first
explain the notion of path complexes (combinatorial generalizations of
facet paths). It will turn out that in the cases we are interested in,
there are many fewer relevant types of path complexes than of simplicial
polytopes.  We then outline the method that allows us to transform
each of the resulting (non)-realizability problems into an (un)-satisfiability
problem. The generation of a particular set of constraints -- the
forbidden shortcuts constraints -- is derived in a dedicated
section. We conclude with a discussion of possible extensions of our
approach. 

\section{Path complexes}
\label{sec:pathcomp}

To restrict the number of cases we need to deal with, we take a closer
look at the facet paths of simplicial polytopes. 
Following the approach
of Bremner, Holt, and Klee~\cite{bremner.ea:pathcomplex}, we consider
the more general setting of pure simplicial complexes of dimension
$d-1$.  The $0$, $d-1$ and $d-2$ simplices are called \emph{vertices},
\emph{facets}, and \emph{ridges}.  We are particularly interested in
the \emph{path complexes} where the dual graph (with \emph{facets} as
nodes and \emph{ridges} shared by two facets as edges) is a path.
It is known~\cites{MR0206823} that when $n\geq 2d$ the maximum
diameter of an $(n,d)$-polytope is always realized by some
\emph{end-disjoint} path (i.e.\ the end vertices do not share a
facet). We thus assume that the start and end facet of all path
complexes considered here are vertex disjoint; we sometimes call such
complexes \emph{end-disjoint path complexes} to emphasize this
feature.

We fix notation for path complexes as follows.  Let $F_0 \dots F_k$ be
the facets. From the definition of a path complex, we know that $F_j$,
$j>0$ can be obtained from $F_{j-1}$ by a \emph{pivot} $(l_j,e_j)$
where $F_j=F_{j-1}\setminus\set{l_j}\cup \set{e_j}$.  Write $F_0 \dots
F_k$ as the rows of $d \times (k+1)$ table so that $l_j$ and $e_j$ are
in the same column.  The list of such column indices is called the
\emph{pivot sequence}.  A basic property of pivot sequences is that
successive indices are distinct
(cf.~\cite{bremner.ea:pathcomplex}*{\S2.1}).

A \emph{non-revisiting} path complex is one where every entering
vertex $e_j$ is distinct; equivalently $|\cup_{j=0}^k F_j|=d+k$.  A
non-revisiting path complex is encoded (up to relabelling of vertices)
in its pivot sequence.  For the general case, we also need to keep
track of \emph{revisits}, i.e.\ when a vertex leaves $F_p$, and
re-enters in $F_q$. 

There are at least two different kinds of symmetry of a path complex
to be considered.  In the first case, we may relabel the vertices of
the initial simplex in $d!$ ways. This symmetry can be removed by
insisting on a particular labelling.  We call a pivot sequence where
the column indices occur in order, i.e.\ the vertex of $F_0$ in
column $c$ leaves before the vertex in column $c+1$, a \emph{canonical
  pivot sequence}.  The second kind of symmetry to be considered is
the choice of initial facet.  In general the two choices may lead to
different canonical pivot sequences. For non-revisiting paths, we
simply keep the lexicographically smaller canonical pivot sequence.

We will next sketch the enumeration of all paths with at most one
revisit as derived by Bremner~et~al.~\cite{bremner.ea:pathcomplex}.
Their first result concerns the number of directed non-revisiting
$d$-paths of length $n$. They show that this number can be expressed
as a Stirling number of the second kind -- we denote these numbers
with $\stirling{n}{k}$.  The basic recursion of the Stirling numbers
can then be used to give a recursive algorithm for generating these paths. The
proof uses one intermediate structure -- \emph{restricted growth
  strings}, i.e.\ $k$-ary strings where the symbols occur in order, and
all $k$ symbols occur. Or put more formally: a restricted growth
string is a sequence $e_{1},\dotsc,e_n$ of symbols from
$\{1,\dotsc,k\}$ such that $e_1=1$ and $e_j=l$ if and only if there
exists an element $e_i$ with $i<j$ such that $e_i=l-1$.

\begin{lemma}[Bremner et al. \cite{bremner.ea:pathcomplex}]
  The number of directed non-revisiting $d$-paths of length $n$ is
  $\stirling{d-1}{n-1}$.
\end{lemma}
\begin{proof}[Sketch of the proof]
  We first argue that there is a bijection between these directed
  non-revisiting $d$-paths of length $n$ and restricted growth strings
  on $n-1$ symbols of length $d-1$. The bijection between facet paths
  and pivot sequences was discussed above. Here we consider the
  correspondence between pivot sequences and restricted growth
  strings.  Given a canonical pivot sequence $p_1 \dots p_n$, we can
  output a restricted growth string $s_1 \dots s_{n-1}$ by setting
  $s_1=1$ and $s_{j-1}$ to the rank of $p_j$ in $\set{1 \dots d}
  \setminus p_{j-1}$ for $j>2$. This transformation is evidently a
  bijection.

  The Stirling numbers of the second kind $\stirling{n}{k}$ count the
  number of partitions of an $n$-element set into $k$ parts. The
  bijection between these partitions and restricted growth strings of
  length $n$ on $k$ elements can be given as follows. Let
  $\pi=(\pi_i)_{i=1}^k$ be a partition of the $n$-element set with $k$
  parts. We may assume that the $k$ parts are ordered according to
  their minimal elements, i.e.\ $\min\pi_1=1$ and for all
  $i\in\{1,\dotsc,k-1\}$ it holds that $\min \pi_i \leq \min
  \pi_{i+1}$. Then we can construct a restricted growth string
  $(e_i)_{i=1}^n$ of length $n$ on $k$ elements by setting $e_i=j$ if
  $i\in\pi_j$. We omit the verification that this mapping is a
  bijection.
\end{proof}

Single revisit paths are generated from non-revisiting paths on one
more vertex by identifying two vertices. We represent this by
partitioning the pivot sequence into three possible empty parts, the
\emph{prefix}, the \emph{loop}, and the \emph{suffix}. The loop
represents the actual revisit where the first pivot is the vertex in
question leaving the facet and the last pivot in the loop is the
vertex returning to the facet. Permissible identifications are
determined by two conditions, given as Lemma~4.2 and Lemma~4.3 in
\cite{bremner.ea:pathcomplex}.  
\begin{lemma}[Bremner et al.~\cite{bremner.ea:pathcomplex}]
\label{lem:4243}
  Let $P|L|S$ be a pivot sequence of a non-revisiting path. Let us
  identify the first and last element of $L$. Then the following
  conditions are necessary for the resulting complex to be an 
  end-disjoint path complex.
  \begin{enumerate}
  \item The loop $L$ must contain three distinct symbols.
  \item Either the first symbol of $L$ must appear in $P$ or the last
    symbol of $L$ must appear in $S$. 
  \end{enumerate}
\end{lemma}

The first condition prevents the creation of a new
ridge. The second condition makes sure that the vertex that is
identified is not on both the first and last facet. It is clear that
these conditions are necessary, thus we do not omit any valid path
complexes by pruning according to them.

For revisiting paths, rather than keeping the lexicographically
smaller canonical pivot sequence, it seems to be better to use a
different symmetry breaking strategy that uses the second condition of
Lemma~\ref{lem:4243}.

\begin{lemma}[Bremner et al.~\cite{bremner.ea:pathcomplex}]
  \label{lem:late-revisit}
  Every combinatorial type of end-disjoint single revisit path has an
  encoding as pivot sequence without a revisit on the first facet.
\end{lemma}
\begin{proof}
  Consider path complex with a pivot sequence $\pi$ with an
  identification in the first facet.  Since the path complex is
  end-disjoint, $\pi$ must not have an identification in the last
  facet according to the second condition of
  Lemma~\ref{lem:4243}. Thus, the reverse pivot sequence has no
  revisit on the first facet. As both of these sequences describe the
  same combinatorial type we may choose the latter.
\end{proof}

As with Lemma~\ref{lem:4243}, it is clear that the condition of (the
proof of) Lemma~\ref{lem:late-revisit} is necessary, and we lose no
combinatorial types of path complexes in filtering by it. Note that in
general our two symmetry breaking strategies for choice of an initial
facet are incompatible and we must choose one.

The final pruning of path complexes that we currently implement is
based on the following result of Bremner, Holt and
Klee~\cite{bremner.ea:pathcomplex}.
\begin{lemma}
\label{lem:not-uniq}
If $\Delta(d-1,n-1) < l-1$, then it will be impossible to embed
geodesically those $(d,l,r)$-paths in which the symbol $1$ appears
uniquely at the beginning or in which $d$ appears uniquely at the end.
\end{lemma}

Note Lemma~\ref{lem:not-uniq} only eliminates candidate complexes when
we have a sufficiently strong bound on $\Delta(d-1,n-1)$. We can
combine the results of this section to generate a set of possible
pivot sequences that satisfy all of the given criteria. In
Table~\ref{tab:enum612} we have listed all pivot sequences for
possible path complexes of length $7$ for polytopes of dimension $6$
on $12$ vertices. We summarize our discussion in the following
proposition. We give an algorithmic way to prove that the necessary
condition outlined in the proposition is true. This is the topic of
the next section.

\begin{proposition}\label{prop:612}
  If all path complexes in Table~2 cannot be realized as a matroid
  polytope, then $\Delta(6,12)=6$ holds.
\end{proposition}

\begin{table}[t]
  \centering
  \caption{Possible pivot sequences in the $(6,12)$ case}
  \begin{tabular}{rrrrrrr}
    (1,7) & (2,8) & (7,9) & (3,10) & (4,7) & (5,11) & (6,12) \\
    (1,7) & (2,8) & (7,9) & (3,10) & (4,11) & (5,7) & (6,12) \\
    (1,7) & (2,8) & (7,9) & (3,10) & (4,11) & (5,12) & (6,7) \\

    (1,7) & (2,8) & (3,9) & (7,10) & (4,11) & (5,7) & (6,12) \\
    (1,7) & (2,8) & (3,9) & (7,10) & (4,11) & (5,12) & (6,7) \\

    (1,7) & (2,8) & (3,9) & (8,10) & (4,11) & (5,8) & (6,12) \\
    (1,7) & (2,8) & (3,9) & (8,10) & (4,11) & (5,12) & (6,8) \\

    (1,7) & (2,8) & (3,9) & (4,10) & (7,11) & (5,12) & (6,7) \\

    (1,7) & (2,8) & (3,9) & (4,10) & (8,11) & (5,12) & (6,8) \\

    (1,7) & (2,8) & (3,9) & (4,10) & (9,11) & (5,12) & (6,9) \\
  \end{tabular}

\label{tab:enum612}
\end{table}

\section{Generation of oriented matroids using SAT solvers}
\label{sec:omsat}

To show that a given path complex cannot be completed to a simplicial
polytope, we show the stronger statement that it cannot be completed
to a matroid polytope, i.e.\ there exists no oriented matroid with the
given path complex in its boundary. This section is devoted to
defining oriented matroids and explaining how the non-existence of
certain oriented matroids (and consequently the non-existence of
certain point sets) can be shown using SAT solvers. The use of SAT
solvers to generate oriented matroids was first described in
\cite{s.diss} and \cite{s.surf}. The method used there is the basis
for the approach outlined in this section. In our computations we used
the SAT solver Minisat by Eén and Sörensson \cite{minisat}.

Oriented matroids have been used before to treat diameter questions of
polytopes; one reference which is particularly interesting is the
thesis of Schuchert \cite{schuchert:diss} where among other things he
confirms that $\Delta(4,11)=6$ in the special case of neighborly
matroid polytopes.

Oriented matroids are a combinatorial abstraction of point
configurations in $\mathbb{R}^d$.  We will use the chirotope axioms of
oriented matroids in the sequel; for further axiom systems and proofs
of equivalence we refer to Chapter~3 of the monograph by
Björner~et~al.~\cite{OM}.  As we are only dealing with simplicial
polytopes, we may always assume that our oriented matroids are
uniform, i.e.\ $\chi(b)\neq 0$ for all $d+1$-sets $b$. This further
simplifies the axioms so that we need to check the following axioms:

\begin{definition}
  Let $E=\{1,\dotsc,n\}$, $r\in\mathbb{N}$, and $\chi: E^r \rightarrow
  \{-1,+1\}$. We call $\mathcal{M}=(E,\chi)$ a uniform oriented
  matroid of rank $r$, if the following conditions are satisfied:
  \begin{enumerate}[(B1)]
  \item The mapping $\chi$ is alternating.
  \item For all $\sigma\in\binom{n}{r-2}$ and all subsets
    $\{x_1,\dotsc,x_4\}\subseteq E\setminus\sigma$ the following
    holds:
    \begin{multline*}
      \{\chi(\sigma,x_1,x_2)\chi(\sigma,x_3,x_4),-\chi(\sigma,x_1,x_3)\chi(\sigma,x_2,x_4),\\\chi(\sigma,x_1,x_4)\chi(\sigma,x_2,x_3)\}=
      \{-1,+1\}
\end{multline*}
  \end{enumerate}
\end{definition}
These relations can be seen as abstractions of the
Grassman-Pl\"{u}cker relations on determinants~\cite{OM}.

We also need to express the fact that the path complex we are
given is in the boundary of the oriented matroid. The fact that an
ordered $d$-set $F$ is a facet of the matroid polytope can be
expressed by enforcing that $\chi(F,e)$ has the same sign for all
$e\in E\setminus F$. We call an oriented matroid where every element
is contained in a facet a \emph{matroid polytope}. As the chirotope
axioms are invariant under negation, we may always assume that the
sign of one base is positive.  Using the fact that the facet incidence graph of
the path complex is connected we can infer the signs of the other
bases that contain a facet of the path complex.
In the case of $n=2d$ (we consider more general situations below),
end-disjointness implies that every point in the oriented matroid must
be on some facet of the input path complex.
This property implies that we do not need to add additional
constraints enforcing convexity 
; the points of the oriented matroid
are already in convex position. To make sure that the starting path
complex is actually a geodesic path in the boundary complex of the
oriented matroid, we need to forbid ``shortcuts'', i.e.\ shorter paths
that connect the end facets of our starting path complex. To enforce this we
need to make sure that for each such possible shortcut at least one
facet is missing. The enumeration of these shortcuts is the subject of
the next section.

As mentioned above we use the approach of \cites{s.diss,s.surf} to
transform these conditions into an instance of SAT. We introduce
variables $[b]$ for each $r$-tuple $b$.  The interpretation of these
variables is that in a satisfying assignment $[b]$ should be true if
and only if $\chi(b)=+$. It follows from condition (B1) that we will
only need $\binom{n}{d}$ of these variables in the final SAT instance.

It turns out that axiom~(B2) yields $16\binom{n}{r-2}\binom{n-r+2}{4}$
CNF constraints. For details we refer to \cites{s.diss,s.surf}. It
remains to explain how to encode the facet path and the forbidden
shortcuts.

Given a $d$-tuple $F$ and an ordering $x_1,\dotsc,x_{n-d}$ of the set
$X=\{1,\dots,n\}\setminus F$, to enforce that $F$ is a facet we need
to add the following clauses: 
\begin{displaymath}
  \left(\bigwedge_{i=1}^{n-d-1}{[F,x_i]\vee
      \neg [F,x_{i+1}]}\right)\wedge
  \left(\bigwedge_{i=1}^{n-d-1}{\neg [F,x_i]\vee
      [F,x_{i+1}]}\right).
\end{displaymath}

Here we use the convention that for a $(d-1)$-tuple
$F=(f_1,\dotsc,f_{d-1})$, the variable $[F,x]$ denotes the variable
$[(f_1,\dotsc,f_{d-1},x)]$.

In order to enforce that some $F\in \mathcal{F}$ is not a facet, we
first construct constraints implied by all $F\in\mathcal{F}$ being on
the boundary, then negate them.  Let $\tau(Y)= (-1)^k$
where $k$ transpositions are required to sort tuple $Y$.

\begin{lemma}
  \label{lem:sigma-back}
  Let $\mathcal{F}=\set{F_1, F_2, \dots F_m}$ be a path complex on the boundary.
  For $i>1$, let   $e_i=F_{i}\setminus F_{i-1}$, $l_i=F_{i-1}\setminus F_{i}$.
  Let $\sigma_1=1$, and for $i>1$ let  $\sigma_i=\tau(F_{i-1},e_i)\tau(F_i,l_i)\sigma_{i-1}$.
  If $F_i$ and $F_{i-1}$ are both facets, then 
  \begin{equation*}
    \sigma_i\chi(F_i, x) = \sigma_{i-1} \chi(F_{i-1}, y) \qquad x \not \in F_i, y \not \in F_{i-1}
  \end{equation*}
\end{lemma}
\begin{proof}
  Suppose that $F_{i-1}$ and $F_i$ are both facets, and let 
  $x\not\in F_i, y \not\in F_{i-1}$, $T_i=F_i \cup F_{i-1}$.
  \begin{align}
    \chi( F_{i-1}, e_i )\chi(F_{i-1}, y) &= \chi(F_i, l_i)\chi(F_i, x) = 1 \notag\\
    \chi( F_{i-1}, e_i )\chi( F_i, l_i)  & = \tau(F_{i-1},e_i)\tau(F_i,l_i)\chi(T_i)^2\notag\\
    &= \tau(F_{i-1},e_i)\tau(F_i,l_i)\notag\\
    \intertext{It follows that}
    \label{eq:tau-back}
    \tau(F_{i-1},e_i)\tau(F_i,l_i)  \chi(F_i, x ) &= \chi(F_{i-1},y) 
\end{align}

The Lemma now follows from the definition of $\sigma_i$, and  \eqref{eq:tau-back}:
\begin{align*}
    \sigma_i\chi(F_i,x) &= \sigma_{i-1} \tau(F_{i-1},e_i)\tau(F_i,l_i) \chi(F_i,x)\notag\\
    &= \sigma_{i-1} \chi(F_{i-1},y).
  \end{align*}
\end{proof}

From Lemma~\ref{lem:sigma-back}, it follows that to force path complex
$\mathcal{F}=\set{F_1, F_2, F_3, \dots F_m}$ not to be on the
boundary, it suffices that
\begin{equation*}
  \set{ \sigma_j\chi(F_j, x) \mid 1\leq j\leq m, x\notin F_j} = \set{+1, -1}
\end{equation*}
Let $z_i(x)$ be the CNF literal corresponding to $\sigma_i \chi( F_i, x)$. Then the two 
corresponding CNF constraints are 
\begin{equation*}
  \left (\bigvee_{i\in \{2\dots |\mathcal{F}|-1\},x\in \{ 1\dots n\}\setminus F_i} z_i(x)\right) \wedge
  \left (\bigvee_{i\in \{2\dots |\mathcal{F}|-1\},x\in \{ 1\dots n\}\setminus F_i} \neg z_i(x)\right) 
\end{equation*}
In our setting we may assume the first and last element of $\mathcal{F}$ are facets.

We get a large number of these clauses and each of them --- on its own
--- is quite weak (as it contains many literals). However, taken
together they lead to the desired contradiction. We note that as we
are looking for a contradiction it is not necessary to generate all of
these clauses.

\section{Shortcuts}
\label{sec:shortcuts}

In the context of testing a $k$ step path complex $\Delta=F_0 \dots
F_{k}$ for geodesic (non)-realizability, we call any path complex (on
the same set of vertices) from $F_0$ to $F_k$ with less than $k$
pivots a \emph{shortcut}. Our general scheme in establishing that a
given path complex is not geodesically realizable is to find a set of shortcuts
$S$, and show that no matroid polytope (and hence no convex polytope)
can contain $\Delta$, but no element of $S$.

A path $\pi=v_0, v_1 \dots v_k$ in graph $G$ is called
\emph{inclusion-minimal} if no proper subset of $v_0 \dots v_k$ is a
path from $v_0$ to $v_k$.  Every inclusion-minimal path is evidently
simple, and every geodesic (shortest path) is inclusion-minimal.

\begin{lemma}
  \label{lem:inclusion-min}
  The inclusion-minimal $(s,t)$-paths of length $k$ in graph $G$ are
  exactly $\pi, t$ where $\pi$ is an inclusion minimal $(s,t')$-path
  of length $k-1$, 
  $t'$ is a neighbour of $t$, and no other $v \in \pi$ is adjacent
  to $t$.
\end{lemma}
\begin{proof}
  Let $\pi=v_0 \dots v_k$ be an inclusion minimal path of length $k$.
  The path $\pi'=v_0 \dots v_{k-1}$ must be inclusion minimal, or
  $\pi$ would not be either. Similarly if $v_k$ is adjacent to some
  $v_j$, $j<k-1$ then $\pi$ is not inclusion minimal.

  Suppose on the other hand we have an inclusion-minimal path
  $\pi'=v_0 \dots v_{k-1}$ such that $v_k$ is adjacent to $v_{k-1}$,
  but not to any other vertex in $\pi'$. Then $\pi=v_0 \dots v_k$ is a
  path from $v_0$ to $v_k$, and if this path is not inclusion minimal,
  then the shortcut must exist in $\pi'$, which is a contradiction.
\end{proof}

From this Lemma we can derive an algorithm that generates all the
potential shortcuts for the given path complex, taking as the graph $G$ the
so-called \emph{pivot graph} whose nodes are $d$-sets and whose edges
are pivots. In the cases we treated so far, we generated all the
shortcuts. As we noted in the previous section, we do not need to make
sure that we have generated all shortcuts to guarantee the correctness
of the result. We therefore implemented a variant of the oriented
matroid search algorithm that 
finds shortcuts in a current candidate realization in the style of
cutting plane algorithms. This is the algorithm we used to prove
Proposition~\ref{prop:612}. It yields non-realizability for all cases
in Table~\ref{tab:enum612}. Thus, we get the following Proposition,
which, together with Proposition~\ref{prop:612}, proves that
$\Delta(6,12)=6$.

\begin{proposition}\label{prop:non-real612}
  None of the path complexes of Table~\ref{tab:enum612} can be
  realized as part of the boundary complex of a matroid polytope. 
\end{proposition}

\section{The case $\Delta(4,11)$}
\label{sec:case-4-11}

The method outlined in the sections above can also be used to show
that $\Delta(4,11)=6$. The difference is the generation of the path
complexes. As we can restrict ourselves to end-disjoint paths, the
number of revisits is bounded by $3$. The paths with up to one revisit
can be generated as outlined in Section \ref{sec:pathcomp}. Here we
can use the symmetry breaking methods outlined in this section. The
paths with two and three revisits are generated similarly, but this
time we get two resp. three loops. Here we did not use any symmetry
reduction. The number of path complexes generated can be found in
Table \ref{tab:411pc}.

\begin{table}
  \centering
  \caption{Number of path complexes after symmetry reduction for the $(4,11)$ case}
  \label{tab:411pc}
  \begin{tabular}{rr}
   \# revisits & \# path complexes\\
   0 & 35\\
   1 & 185 \\
   2 & 354 \\
   3 & 96\\
  \end{tabular}
\end{table}

We can then use the methods outlined in Sections \ref{sec:omsat} and
\ref{sec:shortcuts} to check that none of these path complexes can be
realized as part of a matroid polytope. One complication is that
revisiting paths need not use all of vertices when $n>2d$. Although it
is relatively straightforward to add constraints to enforce that all
points of the oriented matroid are contained in some facet
(see~\cite{bbg} for details) here we rely on the observation that any
realization of a $k$-path which fails to have all of the points on the
boundary is also a realization for fewer points. Since it is
known~\cite{MR0310768} that $\Delta(4,k)\leq 5$, for $k\leq 10$ we may ignore such a
possibility.

\section{Conclusion}
\label{sec:conclusion}

The actual SAT computations took less than one hour for each of the 10
cases in the $(6,12)$-case on a regular desktop computer. However, the
problems get more difficult for higher parameters. The key parameter
seems to be $n-d$; so far, we were not able to finish the computations
for the $(5,12)$ case in reasonable time, even though they are in a
lower dimension. In these cases the SAT solver can be made to produce
an actual proof of infeasibility. However, the produced proofs are too
unwieldy to be checked manually. 

To be able to finish our computation it is crucial that we do not rely
on an explicit enumeration of all matroid polytopes that attain a
given bound. Schuchert \cite{schuchert:diss} enumerated all neighborly
matroid polytopes of dimension $4$ with $11$ vertices, he found
$6\,492$ which had diameter $6$.

For larger computations the SAT problems become considerably
harder. However, it might be more interesting to study the path
complexes in more detail. It would be potentially useful to give
criteria that further reduce the number of equivalence classes of path
complexes to be considered. However, the number of classes will
probably still grow too fast to deal with cases with much larger
parameters.

\begin{table}
\newlength{\tmpboxsixn}\setlength{\tmpboxsixn}{(\tmpbox-\widthof{\textbf{6}})/2}
  \centering
  \caption{Summary of bounds for $\Delta(d,n)$. The bold entries are from the 
  computations described in this article.}
  \label{tab:new}
  \begin{tabular}{ll@{\qquad}c@{\,}c@{\,}c@{\,}c}
     &&\multicolumn{4}{c}{$n-d$}\\
&    & 4 & 5 & 6 & 7\\
\addlinespace
\multirow{4}{*}{$d$}& 4 & \hspace{\tmpboxfour}4\hspace{\tmpboxfour} & \hspace{\tmpboxfive}5\hspace{\tmpboxfive}  & \hspace{\tmpboxfive}5\hspace{\tmpboxfive}  & \hspace{\tmpboxsixn}\textbf{6}\hspace{\tmpboxsixn} \\
&5 & \foll{4} & 5 & 6 & \textbf{\{7,8\}} \\
&6 & \foll{4} & \foll{5} & \textbf{6} & [7,9]\\
&7 & \foll{4} & \foll{5} & \foll{\textbf{6}} & $[7,10]$\\
&
  \end{tabular}
\end{table}

\end{document}